\documentclass{IEEEtran4PSCC}

\usepackage{xspace,amssymb,epsfig,syntonly}
\usepackage{epsfig,amsmath,color}
\usepackage{xspace,syntonly,empheq}
\usepackage{wrapfig}
\usepackage{url}
\usepackage{bbm}
\usepackage{mathtools}
\usepackage{cite}
\usepackage{graphicx}
\usepackage{caption}
\usepackage{subcaption}
\usepackage{algorithm,algorithmic}
\usepackage{verbatim}
\usepackage{amssymb}
\usepackage{amsthm}

\newcommand{\utwi}[1]{\mbox{\boldmath $#1$}}
\newcommand{\E}{{\mathbb{E}}}

\newcommand{\Prob}{{\textrm{Pr}}}

\newcommand{\diag}{{\textrm{diag}}}

\newcommand{\cN}{{\cal N}}

\newcommand{\cR}{{\cal R}}

\newcommand{\cE}{{\cal E}}

\newcommand{\bc}{{\bf c}}
\newcommand{\ba}{{\bf a}}

\newcommand{\be}{{\bf e}}

\newcommand{\bp}{{\bf p}}
\newcommand{\bq}{{\bf q}}

\newcommand{\bs}{{\bf s}}

\newcommand{\bv}{{\bf v}}
\newcommand{\bi}{{\bf i}}

\newcommand{\by}{{\bf y}}

\newcommand{\bB}{{\bf B}}

\newcommand{\bE}{{\bf E}}

\newcommand{\bJ}{{\bf J}}

\newcommand{\bH}{{\bf H}}

\newcommand{\bR}{{\bf R}}

\newcommand{\bI}{{\bf I}}

\newcommand{\bZ}{{\bf Z}}

\newcommand{\bY}{{\bf Y}}

\newcommand{\balpha}{{\utwi{\alpha}}}

\newcommand{\bdelta}{{\utwi{\delta}}}

\newcommand{\bgamma}{{\utwi{\gamma}}}
\newcommand{\bupsilon}{{\utwi{\upsilon}}}

\newcommand{\bPhi}{{\utwi{\Phi}}}

\newcommand{\bSigma}{{\utwi{\Sigma}}}

\newcommand{\bmu}{{\utwi{\mu}}}

\newcommand{\bGamma}{{\utwi{\Gamma}}}

\newcommand{\sfT}{\textsf{T}}

\DeclarePairedDelimiterX{\norm}[1]{\lVert}{\rVert}{#1}

\usepackage{epstopdf}

\begin{document}

\newtheorem{definition}{Definition}
\newtheorem{remark}{Remark}
\newtheorem{proposition}{Proposition}
\newtheorem{lemma}{Lemma}
\def\HS{\hspace{\fontdimen2\font}}
\font\myfont=cmr12 at 16pt

\IEEEoverridecommandlockouts

\title{Efficient Relaxations for Joint Chance Constrained AC Optimal Power Flow  \vspace{-.2cm}} 

\author{Kyri Baker, \emph{Member, IEEE}, and Bridget Toomey}

\maketitle

\begin{abstract}
Evolving power systems with increasing levels of stochasticity call for a need to solve optimal power flow problems with large quantities of random variables. Weather forecasts, electricity prices, and shifting load patterns introduce higher levels of uncertainty and can yield optimization problems that are difficult to solve in an efficient manner. Efficient solution methods for single chance constraints in optimal power flow problems have been considered in the literature; however, joint chance constraints have predominantly been solved via scenario-based approaches or by utilizing the overly conservative Boole's inequality as an upper bound. In this paper, joint chance constraints are used to solve an AC optimal power flow problem which maintain desired levels of voltage magnitude in distribution grids under high penetrations of photovoltaic systems. A tighter version of Boole's inequality is derived and used to provide a new upper bound on the joint chance constraint, and simulation results are shown demonstrating the benefit of the proposed upper bound. 
\end{abstract}


\section{Introduction}
\label{sec:Introduction}
Increasing penetrations of intermittent energy sources in the electric power grid, evolving faster than the corresponding infrastructure, can increase the probability that line congestions may occur, and that voltages may lay outside of desired limits. Rather than considering this randomness as a deterministic input or representing the uncertainty via computationally prohibitive scenario-based approaches, we solve a joint chance constraint problem which prevents overvoltages in the grid with a certain probability. Single chance constraints have been considered in a variety of power systems applications, many addressing the problems of line congestions \cite{Roald13, Bienstock14, Zhang_CC2011}, voltage regulation \cite{BakerNAPS16, DallAneseCDC16}, and energy storage sizing \cite{BakerTSE16}. These works constrain the individual probability of overvoltage at each node, congestion at each line, or bounding individual battery state of charge levels independently; however, perhaps a more relevant constraint to consider is restricting the probability of \emph{all} voltages, line flows, states of charge, etc. being within prescribed limits.

Joint chance constraints have been considered in \cite{Vrak12} for the N-1 security problem, and the joint chance constraint problem was solved by utilizing a sample-based scenario approach. In \cite{Hojjat15}, both single and joint chance constraints were considered to mitigate line congestions, and a Monte-Carlo based approach was developed to estimate the joint probability. Boole's inequality \cite{Boole1854} is a popular choice to provide an upper bound on the original chance constraint \cite{Grosso14, Nemirovski07}, separating the joint chance constraint $P(g_1(x, \delta) \leq 0, ... , g_n(x, \delta) \leq 0) \leq 1 - \epsilon$ into single chance constraints $P(g_1(x, \delta) \leq 0) \leq 1- \epsilon_1 ... P(g_n(x, \delta) \leq 0) \leq 1- \epsilon_n$ for $i = 1 ... n$, where $\sum^{n}_{i = 1} \epsilon_i \leq \epsilon$. A common choice for $\epsilon_i$ is usually $\frac{\epsilon}{n}$ \cite{Grosso14, Nemirovski07}; however, this parameter can also be optimized \cite{Blackmore09}. In \cite{Hong_JCC}, the use of Boole's inequality is avoided by using a Monte Carlo method to solve a sequence of convex optimization problems and compute the joint chance constraint directly; however, it is computationally slow and can only handle relatively small problems.

\begin{figure}[t]
\begin{center}
\includegraphics[width=8.5cm ]{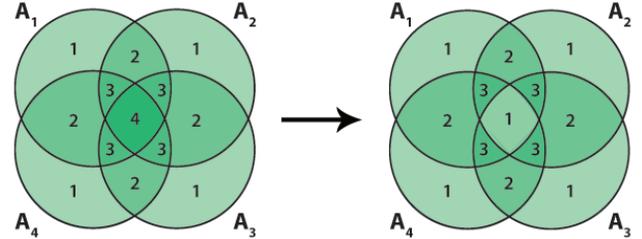}
\end{center}
\vspace{-.3cm}
\caption{Illustrative example: Boole's inequality (left) tends to overestimate the intersection of events; an improved upper bound on $P(A_1 \cup A_2 \cup A_3 \cup A_4)$ is sought by subtracting out the intersection of all events (right).}
\label{Fig:Boole_total}
\vspace{-.2cm}
\end{figure}

In this paper, we provide an improved Boole's inequality which allows for the consideration of a series of single chance constraints, but has the possibility of reducing the conservativeness of the bound provided by Boole. In addition, by exploiting the structure of the voltage regulation problem; i.e., assuming the probability that the system is operating within normal voltage regions is high, the improved bound decreases the cost of the otherwise overly conservative nature of using Boole-based inequalities. It is shown that the new upper bound is tighter than or equal to Boole's inequality, and intuitively amounts to using a bound on the excess probabilities of the intersection of all events using Fr\'{e}chet's inequality, or estimating this intersection with a small number of samples, which Boole's inequality overestimates (see Figure \ref{Fig:Boole_total} for an illustrative example with four events). 

Finally, by utilizing a linearization of the AC power flows in a distribution network, the chance-constrained voltage regulation problem is solved under a high penetration of photovoltaic systems. Simulation results are performed using a modified IEEE-37 node test feeder with community PV systems, and the new bound is compared with a deterministic formulation and the traditional Boole's inequality used in joint chance constraint reformulations. Directions for future work and limitations of the proposed method are also discussed.

\section{Joint Chance Constraint Relaxation}

The joint chance constraint considered here requires the probability of all voltages in the system to be under than the maximum voltage limit with a probability greater than or equal to $1-\epsilon$ :

\begin{equation}
P(g_1(x, \delta) < 0, ... , g_n(x, \delta) < 0) \geq 1 - \epsilon \label{JCC}
\end{equation}

\noindent where $g_1(x, \delta) < 0, ... , g_n(x, \delta) < 0$ constrain the voltage magnitude at each bus $i$, $V_i$, to be less than to the maximum voltage magnitude $\bar{V}$, $x$ is a vector of decision variables, $\delta$ is a jointly distributed Gaussian random vector with mean $\mu$ and positive definite covariance matrix $\Sigma$, and $\epsilon \in (0, 0.5]$. 

\subsection{Improved Boole's Inequality}

Considering each constraint $i$ as an event $B_i$, the joint chance constraint can be written as the intersection of events $P(B_1 \cap B_2 \cap ... \cap B_n)$. Using complementarity, $P(B_1 \cap B_2 \cap ... \cap B_n) = 1- P(B^c_1 \cup B^c_2 \cup ... \cup B^c_n)$. For brevity, define event $B^c_i$ as $A_i$ for each $i$. Boole's inequality states that 
\begin{equation} \label{Booles}
P\Big(\bigcup_{i=1}^{n} A_i\Big) \leq P(A_1) + P(A_2) + ... + P(A_n)
\end{equation}

\noindent for events $A_i$, where $i = 1 ... n$. However, it is clear that the sum of individual probabilities is a conservative upper bound for the union; for example, as seen in Figure \ref{Fig:Boole_total}, the intersection of events is needlessly accounted for multiple times, making the bound defined by Boole's a conservative one. Specifically with regards to the voltage regulation problem, $P(A_1), P(A_2), ... P(A_n)$ refer to the individual probabilities of the voltage at bus $i = 1 ... n$ being greater than or equal to $\bar{V}$.

In this paper, we improve the above inequality by providing a new bound that is equal to or tighter than (\ref{Booles}). First, consider the following:

\begin{align} 
P\Big(\bigcup_{i=1}^{n} A_i\Big) \leq P(A_1) &+ P(A_2) + ... + P(A_n)\nonumber \\
& - P(A_1 \cap A_2 \cap ... \cap A_n) \label{New_1}
\end{align}
\noindent Which subtracts the intersection of all events. It is clear that $P(A_1 \cap A_2 \cap ... \cap A_n) \geq 0$, and subtracting this from the sum of the individual probabilities does not remove any portion of the feasible region of the original union, which will be shown later. We propose (\ref{New_1}) as an equal or tighter bound on the union than (\ref{Booles}). Considering the fact that if this union is nonzero, it will be repeated $n$ times, with $n-1$ of those instances being unnecessary. Thus, we can extend (\ref{New_1}) to the following:

\begin{align} 
P\Big(\bigcup_{i=1}^{n} A_i\Big) \leq P(A_1) &+ P(A_2) + ... + P(A_n)\nonumber \\
& - (n-1)P(A_1 \cap A_2 \cap ... \cap A_n) \label{New_2}
\end{align}

\noindent If this intersection cannot be computed, by using Fr\'{e}chet's inequality, we can provide an alternate upper bound:

\begin{align} 
P\Big(\bigcup_{i=1}^{n} A_i\Big) \leq \hspace{.1cm} &P(A_1) + P(A_2) + ... + P(A_n) - \nonumber \\
&(n-1)[P(A_1) + ... + P(A_n) - (n-1)]_+ \label{New_Booles}
\end{align}

\noindent Where $[\cdot]_+$ denotes $max(0, \cdot)$. \\

\noindent \textbf{Proposition 1.} \emph{The bounds provided by (\ref{New_2}) and (\ref{New_Booles}) provide a valid upper bound that is equal to or tighter than (\ref{Booles}).} \\

\begin{proof} 
First, we prove by induction that the following inequality holds for any $n = 2$ (for $n=1$, the joint chance constraint becomes a single chance constraint):

\begin{small}
\begin{align}
P(A_1 \cup A_2) \leq P(A_1) + P(A_2) - (n-1)P(A_1 \cap A_2) 
\end{align}
\end{small}

\noindent By definition, $P(A_1 \cup A_2) = P(A_1) + P(A_2) - P(A_1 \cap A_2)$; thus, the base case holds. Now, assume that case $n$ holds true:

\begin{small}
\begin{align*}
P\Big(\bigcup_{i=1}^{n} A_i\Big) \leq \sum_{i=1}^{n}P(A_i) - (n-1)P\Big(\bigcap_{i=1}^{n}A_i\Big)
\end{align*}
\end{small}

\noindent we wish to prove that the inequality still holds for $n+1$; that is

\begin{small}
\begin{align*}
P\Big(\bigcup_{i=1}^{n+1} A_i\Big) \leq \sum_{i=1}^{n+1}P(A_i) - nP\Big(\bigcap_{i=1}^{n+1}A_i\Big)
\end{align*}
\end{small}

\noindent By the base case,

\begin{small}
\begin{align}
P\Big(\bigcup_{i=1}^{n+1} A_i\Big) \leq \Big(P\Big(\bigcup_{i=1}^n A_i\Big) + P(A_{n+1})\Big) - P\Big(\bigcap_{i=1}^{n}A_i \cap A_{n+1}\Big) \label{basecase}
\end{align}
\end{small}

\noindent But

\begin{small}
\begin{align*}
\bigcap_{i=1}^{n+1}A_i \subseteq  \Big(\bigcap_{i=1}^{n}A_i \Big) \cap A_{n+1}
\end{align*}
\end{small}

\noindent and thus

\begin{small}
\begin{align*}
P\Big(\bigcap_{i=1}^{n+1}A_i\Big) \leq  \Big(P\Big(\bigcap_{i=1}^{n}A_i \Big) \cap A_{n+1}\Big)
\end{align*}
\end{small}

\noindent which, by the inductive hypothesis, yields the upper bound 

\begin{small}
\begin{align*}
P\Big(\bigcap_{i=1}^{n+1}A_i\Big) &\leq P\Big(\bigcap_{i=1}^n A_i\Big) + P(A_{n+1}) - P\Big(\bigcap_{i=1}^{n+1}A_i\Big)\\
& \leq \sum_{i=1}^nP(A_i) - (n-1)P\Big(\bigcap_{i=1}^n A_i \Big) \\
& \hspace{1.8cm}+ P(A_{n+1}) - P\Big(\bigcap_{i=1}^{n+1} A_i\Big)
\end{align*}
\end{small}

\noindent and because $\cap_{i=1}^{n+1}A_i \subseteq \cap_{i=1}^n A_i$, $P(\cap_{i=1}^{n+1}) \leq P(\cap_{i=1}^n A_i)$. Therefore, we can write

\begin{small}
\begin{align*}
P\Big(\bigcap_{i=1}^{n+1}A_i\Big) &\leq \sum_{i=1}^n P(A_i) + P(A_{n+1}) \\
&\hspace{1cm} - (n-1)P\Big(\bigcap_{i=1}^{n+1}A_i\Big) - P\Big(\bigcap_{i=1}^{n+1}A_i\Big)\\
& \leq \sum_{i=1}^{n+1}P(A_i) - (n)P\Big(\bigcap_{i=1}^{n+1}A_i\Big)
\end{align*}
\end{small}

\noindent Lastly, by Fre\'{e}chet's inequality, we have

\begin{small}
\begin{align*}
P\Big(\bigcap_{i=1}^{n+1}A_i\Big) &\leq \sum_{i=1}^n P(A_i) + n\Big[\sum_{i=1}^{n+1}P(A_i) - n \Big]_+
\end{align*}
\end{small}

\noindent and the proof by induction is complete. For (\ref{New_2}) and (\ref{New_Booles}) to be equal or tighter upper bounds than (\ref{Booles}), it is sufficient to show that 

\begin{small}
\begin{align}
\sum_{i=1}^{n}P(A_i) - (n-1)\Big[\sum_{i=1}^{n}P(A_i)  - (n-1)\Big]_+ \leq \sum_{i=1}^{n}P(A_i) \label{lem1}
\end{align}
\end{small}

\noindent Because $n>1$ and $[\sum_{i=1}^{n}P(A_i)  - (n-1)]_+ \geq 0$, it is clear that the left hand side of (\ref{lem1}) is always less than or equal to the right hand side. For small $\sum_{i=1}^{n}P(A_i)$ or large $n$, this bound will likely end up being equivalent to Boole's inequality, because $[\sum_{i=1}^{n}P(A_i)  - (n-1)]_+$ is likely to be $0$. However, for sensitivity studies, such as identifying under which situations the probability of overvoltage is high, this bound could prove useful. For the application considered in this paper, we will utilize the bound (\ref{New_2}) which directly considers the intersection.

\end{proof}

\subsection{Reformulating the Joint Chance Constraints} \label{sec:JCC}

Consider the original joint chance constraint and its complement:

\begin{eqnarray*}
& \hspace{-1cm}P(g_1(x, \delta) < 0 \cap ... \cap g_n(x, \delta) < 0) \iff  \\
& \hspace{.4cm} 1 - P(g_1(x, \delta) \geq 0 \cup ... \cup g_n(x, \delta) \geq 0)
\end{eqnarray*}

\noindent The probability of the union of events can thus be written as

\begin{eqnarray*}
P(g_1(x, \delta) \geq 0 \cup ... \cup g_n(x, \delta) \geq 0) \leq \epsilon
\end{eqnarray*}

\noindent Then, the final joint chance constraint reformulation can be written as the following series of single chance constraints:

\begin{eqnarray*}
&P(g_1(x,&\hspace{-.25cm}\delta) \geq 0) \leq \epsilon_1 \\
&P(g_2(x,&\hspace{-.25cm}\delta) \geq 0) \leq \epsilon_2 \\
&&\vdots \\
&P(g_n(x,&\hspace{-.25cm}\delta) \geq 0) \leq \epsilon_n
\end{eqnarray*}

\noindent where according to (\ref{New_2}),

\begin{eqnarray} \label{ep_bound}
\sum_{i=1}^n \epsilon_i - P\Big(A_1 \cap ... \cap A_n\Big) \cdot (n-1) \leq \epsilon.
\end{eqnarray}

\noindent  It will be shown in the following section that the intersection $P\Big(A_1 \cap ... \cap A_n\Big)$ can be efficiently estimated using a small number of samples.

\section{System model and Linearization}
\label{sec:models_and_approx}

\subsection{Distribution Network}
\label{sec:systemmodel}

Consider a distribution feeder comprising of $N+1$ nodes collected in the
set $\cN \cup \{0\}$, $\cN := \{1,\ldots,N\}$, and lines represented by the set of
edges $\cE := \{(m,n)\} \subset \cN \times \cN$. Let $V_n \in \mathbb{C}$ and $I_n \in \mathbb{C}$ denote the phasors for the line-to-ground voltage and the current injected at node $n \in
\cN$, respectively, and define the $N$-dimensional complex vectors  $\bv := [V_1, \ldots, V_N]^\sfT \in \mathbb{C}^{N}$ and $\bi := [I_1, \ldots, I_N]^\sfT \in
\mathbb{C}^{N}$. On the other hand, node $0$ denotes the secondary
of the distribution transformer, and it is taken to be the slack bus. Using Ohm's and Kirchhoff's circuit laws, the following linear relationship can be established:
\begin{align}
\label{eq:iyv}
\left[
\begin{array}{c}
I_0 \\
\bi
\end{array}
\right] = 
\underbrace{\left[
\begin{array}{cc}
 y_{00}  & \bar{\by}^\sfT  \\
 \bar{\by}  & \bY
\end{array}
\right]}_{:= \bY_{\mathrm{net}}}
\left[
\begin{array}{c}
V_0 \\
\bv
\end{array}
\right]
\end{align}
where the system admittance matrix $\bY_{\mathrm{net}} \in \mathbb{C}^{N+1 \times N+1}$ is formed
based on the system topology and the $\pi$-equivalent circuit of the distribution lines (see e.g.,~\cite[Chapter 6]{kerstingbook} for additional details on distribution line modeling), and is partitioned in sub-matrices with  the following dimensions: $\bY \in \mathbb{C}^{N \times N}$, $\overline \by \in \mathbb{C}^{N \times 1}$, and $y_{00} \in \mathbb{C}$. The voltage at the slack bus is defined as $V_0 = \rho_0 e^{\mathrm{j} \theta_0}$, with $\rho_0$ denoting the voltage magnitude at the secondary of the step-down transformer. Lastly, $P_{\ell,n}$ and $Q_{\ell,n}$ denote the real and reactive demands at node $n \in \cN$, and define the vectors  $\bp_\ell := [P_{\ell,1}, \ldots, P_{\ell,N}]^\sfT$ and $\bq_\ell := [Q_{\ell,1}, \ldots, Q_{\ell,N}]^\sfT$; if no load is present at node $n \in \cN$, then $P_{\ell,n} = Q_{\ell,n} = 0$, $\forall \, t$. 

\subsection{PV Systems}
Random quantity $P_{\textrm{av},n}$ denotes the maximum renewable-based generation at node $n \in \cN_R \subseteq \cN$  -- hereafter referred to as the available real power. Particularly, $P_{\textrm{av},n}$ coincide with the maximum power point at the AC side of the inverter. When RESs operate at unity power factor and inject the available real power $P_{\textrm{av},n}$,  a set of challenges related to power quality and reliability in distribution systems may emerge for sufficiently high levels of deployed RES capacity~\cite{Liu08}. For example, overvoltage at a particular node may be experienced when RES generation exceeds the load of that consumer~\cite{Liu08}.  Efforts to ensure reliable operation of existing distribution systems with increased behind-the-meter renewable generation are focus on the possibility of inverters providing reactive power compensation and/or curtailing real power. To account for the ability of the RES inverters to adjust the output of real power, let $\alpha_n \in [0,1]$ denote the fraction of available real power curtailed by RES-inverter $n$. If no PV system/inverter is at a particular node $i$, $P_{av,i} = \alpha_i = 0$. For convenience, define the vectors $\balpha :=[\alpha_1,\ldots, \alpha_N  ]^\sfT$ and $\bp_{\textrm{av}} := [P_{\textrm{av},1}, \ldots, P_{\textrm{av},N}]^\sfT$where $\alpha_n = 0$, $P_{\textrm{av},n} = 0$, and ${Q}_{n} = 0$ for $n \in \cN \backslash \cN_R$.

The available real power from solar is modeled as $\bp_{\textrm{av}} = \bar{\bp}_{\textrm{av}} + \bdelta_\textrm{av}$, where $\bar{\bp}_{\textrm{av}} \in \mathbb{R}^N$ is a vector of the forecasted values and $\bdelta_\textrm{av}  \in \cR_\textrm{av}  \subseteq \mathbb{R}^{N}$ is a random vector whose distribution captures spatial dependencies among forecasting errors. We assume that the distribution system operator has a certain amount of information about the probability distributions of the forecasting errors $\bdelta_\textrm{av}$. This information can come in the form of either knowledge of the  probability density functions, or a model of $\bdelta_\textrm{av}$ from which one can draw samples. In this paper, we make the assumption that these errors are Normally distributed; however, distributionally robust formulations of single chance constraints \cite{Tyler15, Nemirovski07} can easily be incorporated into the framework here.

\subsection{AC Power Flow Approximation}
\label{sec:approximate}

Using~\eqref{eq:iyv}, the net complex-power injections can be compactly written as
\begin{align}
\bs = \mathrm{diag}\left(\bv \right) \left(\bY^* (\bv)^* + \overline{\by}^* (V_0)^* \right).
\label{eq:balance}
\end{align}
where $\bs := [\bs_1, \ldots, \bs_N]^\sfT$ and $S_i = (1 - \alpha_i) P_{\textrm{av},i} - P_{\ell,i}  - \textrm{j} (Q_{\ell,i})$. This equation typically appears in the form of a constraint in standard formulations of the OPF task, and renders the underlying optimization problem nonconvex~\cite{LavaeiLow}. Non-convexity implies that off-the-shelf solvers for nonlinear programs may not achieve global optimality; from a computational standpoint, their  complexity may become prohibitive with the increasing of the problem size~\cite{PaudyalyISGT}. Semidefinite relaxation techniques have been employed to bypass the nonconvexity of voltage-regulation and power-balance constraints, and yet achieve globally optimal solutions of the nonconvex OPF under a variety of conditions (see e.g.,~\cite{LavaeiLow}). Several other convex relaxation techniques have also been investigated (see e.g.,~\cite{zhang2013relaxed,madani2014,coffrin2015qc} and pertinent references therein) many of which could also be utilized in this same context. Here, to derive a convex reformulation of the chance constrained OPF, linear surrogates of~\eqref{eq:balance} and voltage-regulation constraints will be utilized next.

To this end, collect voltages $\{V_n\}_{n \in \cN}$ in the vector $\bv := [V_1, \ldots, V_N]^\sfT \in \mathbb{C}^{N}$ and the voltage magnitudes $\{|V_n|\}_{n \in \cN}$ in $\bv_{\ba} :=  [|V_1|, \ldots, |V_N|]^\sfT \in \mathbb{R}^{N}$. The objective is to obtain  approximate power-flow relations whereby voltages  are \emph{linearly} related to injected powers as 
\begin{align} 
\bv & \approx \bH \bp + \bJ \bq + \bc \label{eq:approximateV} \\
\bv_{\ba} & \approx \bR \bp + \bB \bq + \ba, \label{eq:approximate} 
\end{align}
where $\bp := \Re\{\bs\}$ and $\bq := \Im\{\bs\}$~\cite{sairaj2015linear,swaroop2015linear}. By using this relationship,  voltage constraints $|V_i| \leq V^{\mathrm{max}}$, $i \in \cN$,  can be approximated as $ \bR \bp + \bB \bq + \ba \preceq V^{\mathrm{max}} \mathbf{1}_N$, while~\eqref{eq:approximateV}-\eqref{eq:approximate} represents surrogates of~\eqref{eq:balance}. Following~\cite{sairaj2015linear,swaroop2015linear}, the matrices $\bR, \bB, \bH, \bJ$ and the vectors $\ba, \bc$ are obtained as follows.  

Consider linearizing the AC power-flow equation around a given voltage profile $\bar{\bv} := [\bar{V}_1, \ldots, \bar{V}_N]^\sfT$~\cite{sairaj2015linear,swaroop2015linear}. In the following,  the voltages $\bv$ satisfying the nonlinear power-balance equations~\eqref{eq:balance} are expressed as $\bv = \bar{\bv} + \be$, where the entries of $\be$ capture deviations around the linearization points $\bar{\bv}$. Define  $\bar{\bv}_a \in \mathbb{R}^N_+$ the magnitudes of voltages $\bar{\bv}$, and let $\bar{\bgamma} \in \mathbb{R}^N$ and $\bar{\bmu} \in \mathbb{R}^N$ collect elements $\{\cos( \bar{\theta}_i)\}$ and $\{\sin( \bar{\theta}_i)\}$, respectively,  where $\bar{\theta}_i$ is the angle of the nominal voltage $\bar{V}_i$.  Expanding on~\eqref{eq:balance}, and discarding  second-order terms such as $\mathrm{diag}\left(\be\right) \bY^* \be^*$, it turns out that~\eqref{eq:balance} can be approximated as $\bGamma  \be + \bPhi \be^* = \bs + \bupsilon$, where $\bGamma := \mathrm{diag}\left(\bY^* \bar{\bv}^* + \overline{\by}^* V_0^* \right)$, $\bPhi := \mathrm{diag}\left(\bar{\bv}\right) \bY^*$,  and $\bupsilon := -\mathrm{diag}\left(\bar{\bv}\right) \left(\bY^* \bar{\bv}^*  + \overline{\by}^* V_0^*\right)$. Next, consider then the following choice of the nominal voltage $\bar{\bv}$:
\begin{equation} 
\label{eq:V_noload}
\bar{\bv} = - \bY^{-1} \overline{\by} V_0 \, .
\end{equation}
Using~\eqref{eq:V_noload}, it follows that $\bGamma = \mathbf{0}_{N \times N}$ and $\bupsilon = \mathbf{0}_N$, and therefore one  obtains the linearized power-flow expression 
\begin{equation} \label{eq:linearized-power-flow}
\mathrm{diag}\left(\bar{\bv}^*\right) \bY \be = \bs^*. 
\end{equation}
Notice that  matrix $\bY$ is diagonally dominant and irreducible~\cite{sairaj2015linear}. Particularly, it is diagonally dominant by construction since $|y_{kk}| \geq \sum_{i \neq k} |y_{ki}|$ for all $i \in \cN$; it is also irreducibly diagonally dominant if $|y_{0k}| > 0$ for any $k$. Then, a solution to~\eqref{eq:linearized-power-flow} can be expressed as $\be =  \bY^{-1}\mathrm{diag}^{-1}(\bar{\bv}^*) \bs^*$. Thus, expanding on this relation,  the approximate voltage-power relationship~\eqref{eq:approximateV} can be obtained by defining the matrices:   
\begin{subequations} 
\label{eq:Param_approx}
\begin{align}
& \hspace{-.2cm}  \bar{\bR}  = \bZ_R \diag(\bar{\bgamma}) (\diag(\bar{\bv}_{\ba}))^{-1} -\bZ_I \diag(\bar{\bmu}) (\diag(\bar{\bv}_{\ba}))^{-1} \\
&  \hspace{-.2cm} \bar{\bB}  = \bZ_I \diag(\bar{\bgamma}) (\diag(\bar{\bv}_{\ba}))^{-1} + \bZ_R \diag(\bar{\bmu}) (\diag(\bar{\bv}_{\ba}))^{-1}  
\end{align}
\end{subequations}
where $\bZ_R := \Re\{\bY^{-1}\}$ and $\bZ_I := \Im\{\bY^{-1}\}$, and setting $\bH = \bar{\bR} + \mathrm{j} \bar{\bB}$, $\bJ = \bar{\bB} - \mathrm{j} \bar{\bR}$, and $\bc = \bar{\bv}$ . If the entries of $\bar{\bv}$ dominate those in $\be$, then $\bar{\bv}_{\ba} + \Re\{\be\}$ serves as a first-order approximation to the voltage magnitudes across the distribution network~\cite{sairaj2015linear}, and  relationship~\eqref{eq:approximate} can be obtained by setting $\bR = \bar{\bR}$, $\bB = \bar{\bB}$, and $\ba = \bar{\bv}_a$. Equations~\eqref{eq:approximateV}--\eqref{eq:approximate}  are now utilized to solve the relaxed joint chance constraint problem.

\section{Chance Constrained Formulation}
\label{sec:optimization}
\subsection{Optimization problem reformulation}
The original, unrelaxed joint chance constraint optimization for voltage regulation in distribution systems shown below:

\begin{subequations} 
\label{PJCC}
\begin{align} 
 \textrm{(P0)}&  \min_{\substack{\bv_{\ba}, \balpha}}  \hspace{.2cm}  \E(f(\bv_{\ba}, \balpha, \bp_{\ell}, \bq_{\ell})) \label{eq:Pmg_cost}\\
& \mathrm{subject\,to} \nonumber \\
& \bv_{\ba}  = \bR ((\bI - \diag\{\balpha\}) \bp_{\textrm{av}} -  \bp_{\ell}) \\
& \hspace{1cm}- \bB \bq_{\ell} + \ba\label{mg-balance-t} \\
& \Prob\{v_{a,1} \leq V_{\mathrm{max}}, ... , v_{a,n} \leq V_{\mathrm{max}} \} \geq 1- \epsilon \label{mg-prob-jcc}\\
& 0 \leq \alpha_i \leq 1 \hspace{3.05cm}  \label{mg-alpha} 
\end{align}
\end{subequations}
\noindent for all $i \in \cN$, where $v_{a,k}$ denotes the $k$-th element of $\bv_{\ba}$. Constraint~\eqref{mg-balance-t} represents a surrogate for the power balance equation; constraint \eqref{mg-prob-jcc} is the joint chance constraint that require every voltage magnitude in the grid to be within upper and lower limits with at least $1-\epsilon$ probability; and constraint \eqref{mg-alpha} limits the curtailment percentage from $0-100\%$. The cost function $f(\bv_{\ba}, \balpha, \bp_{\ell}, \bq_{\ell})$ is convex and can consider a sum of penalties on curtailment, penalties on power drawn from the substation, penalties on voltage violations, among other objectives\footnote{A more generalized voltage regulation formulation would also include lower limits on voltage; for simplicity of exposition in the high-penetration PV case shown in this paper, only upper limits were considered.}.

As derived in Section \ref{sec:JCC}, the above optimization problem, which is in general nonconvex due to constraint (\ref{mg-prob-jcc}), can be rewritten as the following optimization problem with single chance constraints:

\begin{subequations} 
\label{PCC}
\begin{align} 
 \textrm{(P1)}&  \min_{\substack{\bv_{\ba}, \balpha}}  \hspace{.2cm}  \E(f(\bv_{\ba}, \balpha, \bp_{\ell}, \bq_{\ell})) \\
& \mathrm{subject\,to} \nonumber \\
& \bv_{\ba}  = \bR ((\bI - \diag\{\balpha\}) \bp_{\textrm{av}} -  \bp_{\ell}) \\
& \hspace{1cm}- \bB \bq_{\ell} + \ba \\
& \Prob\{v_{a,i} \leq V_{\mathrm{max}} \} \geq 1- \epsilon_i \label{mg-prob-cc}\\
& 0 \leq \alpha_i \leq 1 \hspace{3.05cm}  
 \end{align}
\end{subequations}

\noindent for all $i \in \cN$, and each $\epsilon_i$ is chosen such that (\ref{ep_bound}) holds.

\subsection{Analytical Reformulation of Single Chance Constraints}

The constraints (\ref{mg-prob-cc}) can be then reformulated as exact, tractable constraints \cite{BoVa04}, assuming $\epsilon \leq 0.5$. Assuming the joint distribution of the random variables is a multivariate Gaussian with  mean $\bmu$ and covariance matrix $\bSigma$, define $\mu_i$ as the $i$th value in $\bmu$ and $\sigma$ as the $(i,i)$ entry in $\bSigma$. Then, define the following function at each node $i \in \cN$:

\begin{align*}
h(p_{av,i}) = &\sum_j(R_{ij}[(1-\alpha_j)p_{av,j} - p_{\ell, j}]) \\
&- \sum_j(B_{ij}q_{l,j}) + a_i - V_{max}
\end{align*}

\noindent where $R_{ij}$ is the $(i,j)$th entry of $\bR$, $B_{ij}$ is the $(i,j)$th entry of $\bB$, and $a_i$ is the $i$th element of $\ba$. Then $h(p_{av,i})$ is also Normally distributed with the following mean $\mu'_i$ and variance $\sigma'_i$:

\begin{align*}
\mu'_i &= \sum_j(R_{ij}[(1-\alpha_j)\mu_j - p_{\ell, j}]) - \sum_j(B_{ij}q_{l,j}) + a_i - V_{max}\\
\sigma'_i &=  \sum_jR_{ij}(1-\alpha_j)\sigma_j
\end{align*}

\noindent Thus, the constraints (\ref{mg-prob-cc}) can be reformulated using the Gaussian cumulative distribution function (CDF) $\Phi$:

\begin{align*}
\Prob\{h(p_{av,i}) \leq 0\} = \Phi\Big(\frac{0-\mu'_i}{\sigma'_i}\Big) \geq 1-\epsilon_i
\end{align*}

\noindent With the final analytical constraint written using the quantile function (the inverse of the Gaussian CDF):

\begin{align}
R_i[(1-\alpha_i)\mu_i - p_{\ell,i}] - B_iq_{\ell,i} + a_i - V_{max} \nonumber \\
\leq -R_i\alpha_i\sigma_i\Phi^{-1}(1-\epsilon_i) \label{final_CC}
\end{align}

\noindent Which can be explicitly included into problem (P1) for each node $i \in \cN$ in place of constraints (\ref{mg-prob-cc}).

\subsection{Estimating the probability of intersection}

Using a Monte-Carlo approach to reformulate the joint chance constraint into a series of deterministic ones can provide large computational burdens. According to \cite{campi}, the number of realizations $N_m$ of the uncertain parameter that need to be included as constraints in place of the chance constraint, with a $\beta$ percent confidence level, should be at least the following:

\begin{align}
N_m \geq \frac{2}{\epsilon}\Big( ln\Big( \frac{1}{\beta} \Big) + N_d \Big)
\end{align}

\noindent where $N_d$ is the total number of decision variables. For example, if $\epsilon = 0.01$, $N_d = 37$, and we desire with a $99\%$ confidence level that the chance constraint will be fulfilled with probability $1-\epsilon$, 8,321 realizations of the random parameter is needed at each timestep. If we were to consider a model predictive control approach to the voltage regulation problem as in our previous work \cite{DallAneseCDC16}, a two-hour horizon with five minute timesteps, with 37 variables per timestep at these same confidence and violation levels would require $N_m \geq 178,521$ deterministic constraints at each receding horizon optimization. Here, we do not propose Monte-Carlo sampling for transforming the chance constraints into thousands of deterministic ones. Instead, we propose to tighten the initial bound used to determine each $\epsilon_i$ by estimating the probability $P(A_1 \cap ... \cap A_n)$ through Monte-Carlo simulations.

\begin{figure}[t]
\begin{center}
\includegraphics[width=7.5cm ]{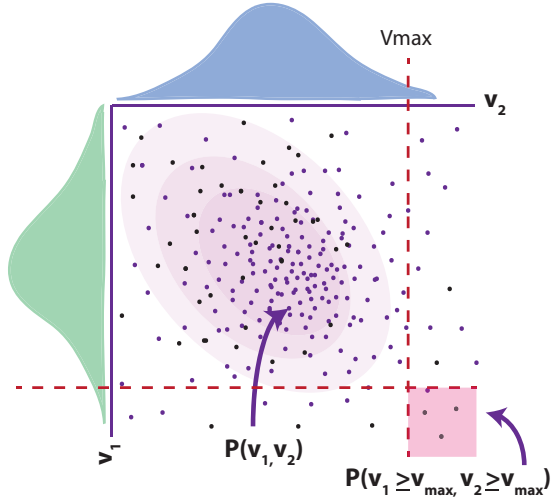}
\end{center}
\vspace{-.3cm}
\caption{Illustrative example with a joint multivariate normal distribution of voltage at two nodes and the region of interest which is the probability of an overvoltage condition at both nodes.}
\label{fig:joint_cc}
\vspace{-.2cm}
\end{figure}

An illustrative example of a joint multivariate normal distribution representing the distribution of voltage magnitude at two nodes is shown in Figure \ref{fig:joint_cc}, illustrating this region of interest. The proposed approach has the potential to be an improvement over Boole's, due to the fact that we can subtract this intersection from the sum of marginal probabilities, but if computational time is important, it is also an improvement over Monte-Carlo methods that estimate the joint chance constraint directly and require a large number of samples in order to accurately represent the original constraint. For each Monte-Carlo sample $m$ generated from the joint multivariate normal distribution, we compute the intersection as follows:

\begin{align} \label{MC}
P(A_1 \cap ... \cap A_n) \approx \frac{\sum_{m=1}^{N_m} \textbf{1}_{\bv_{\ba} \geq V_{max}}(\bv_{\ba}(m))}{N_m}
\end{align}

\noindent That is, for each sample that is drawn from the joint distribution of $\bp_{av}$, the voltage magnitude is evaluated to see if $V_{max}$ is violated. The policy for $\balpha$ is utilized from the solution of that timestep where each $\epsilon_i = \frac{\epsilon}{|\cN_R|}$, which is an upper bound for the new $\epsilon_i$ which will be updated after (\ref{MC}) is performed and calculated using (\ref{ep_bound}). This is a fast calculation; for example, in the simulations in the following section, performing 10,000 calculations only required 0.19 seconds, and the resulting intersection probability calculated with 10,000 samples and with 100,000 samples only differed by .01\%. 
\section{Illustrative Simulations}
\label{sec:sims}

\subsection{Results on IEEE-37 Node Test Feeder}

\begin{figure}[t]
\begin{center}
\includegraphics[width=7cm ]{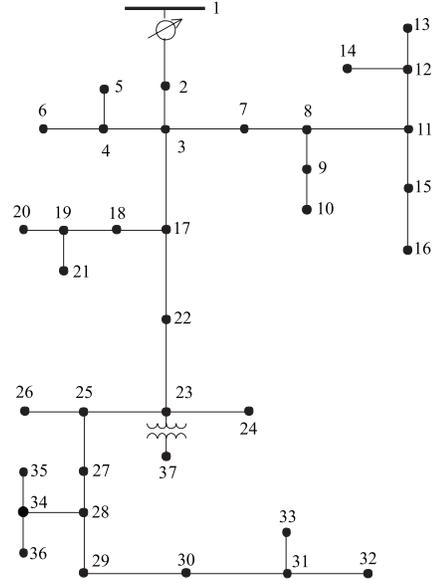}
\end{center}
\vspace{-.3cm}
\caption{The IEEE 37-node distribution feeder.}
\label{fig:37node}
\vspace{-.2cm}
\end{figure}
\begin{figure}[t]
\begin{center}
\includegraphics[width=9.5cm ]{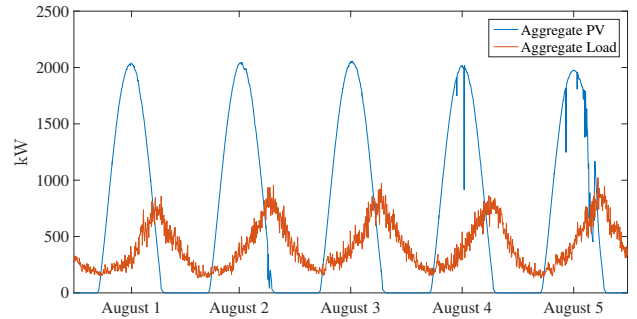}
\end{center}
\vspace{-.3cm}
\caption{Total feeder loading and available PV power during August 1-5, 2012.}
\label{fig:totLoadPV}
\vspace{-.2cm}
\end{figure}

The IEEE-37 node test feeder \cite{IEEE37}, as seen in Figure \ref{fig:37node}, was used for the following simulations. The actual five-minute load and solar irradiance data was obtained from \cite{Bank13} for the simulations, and shown in Figure \ref{fig:totLoadPV}. In order to emulate a situation with high-PV penetration and risks of overvoltage, 16 $200$-kW rated PV systems were placed at nodes 3-18. The considered cost function seeks to minimize renewable curtailment; specifically, 

\begin{align} 
f(\bv_{\ba}, \balpha, \bp_{\ell}, \bq_{\ell}) =  \sum_{i \in \cN} b_i  \alpha_i^2 ,  \label{eq:cost_sim}
\end{align} 

\noindent where the cost of curtailing power at each node is set to be $b_i = \$0.10$. The number of samples was set to $N_m = 100,000$, which, as stated above, can be computed in a fraction of a second. The considered joint chance constraint considers maintaining voltages at nodes $3...37$, with the substation voltage at node $1$ fixed to $1.03$ p.u. Each $\mu_i$, $i=1...n$ was chosen to be the power generated from the forecasted PV at that node, based on the shape of the aggregate solar irradiance from \cite{Bank13} and shifted using samples from a uniform distribution from +/- 1 kW across each node. The covariance matrix $\Sigma$ was formed by setting each entry $(i,j)$ to $\Sigma_{ij} = \bE[(P_{av,i} - \mu_i)(P_{av,j} - \mu_j)^T]$. Non-Gaussian probability distributions can also be considered in the proposed methodology by considering distributionally robust convex approximations of single chance constraints \cite{BakerNAPS16, DallAneseCDC16, Tyler15}.

\subsection{Comparison with Deterministic and Boole's Formulations}
Simulations were performed using data from the first five days of August 2012, and using the joint chance constraint violation parameter $\epsilon = 0.01$.

\subsubsection{Value of intersection term}
\begin{figure}[t!]
\begin{center}
\includegraphics[width=9.5cm ]{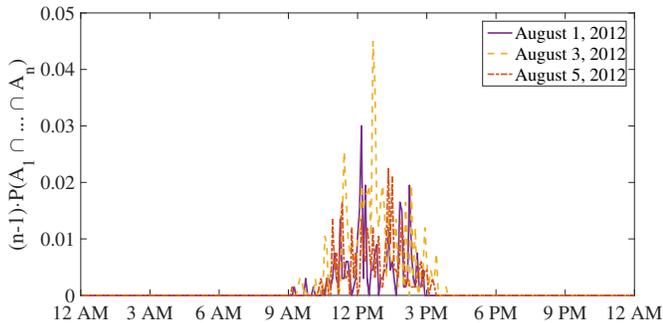}
\end{center}
\vspace{-.3cm}
\caption{Value of intersection term estimated from 100,000 Monte-Carlo samples.}
\label{fig:P_intersection}
\vspace{-.2cm}
\end{figure}

During times of no solar irradiance, the observed probability of overvoltage is low; hence, the term $P(A_1 \cap ... \cap A_n)$ was not observed to be nonzero during these times. The value of this term can be seen Figure \ref{fig:P_intersection} for August 1, 3, and 5; although the probability of intersection is rather small during peak solar irradiance hours, multiplying this term by $(n-1)$ results in a more significant value, reducing the conservativeness of Boole's inequality. Each $\epsilon_i$ in Boole's case was chosen to be $\epsilon_i = \frac{\epsilon}{n}$; and in the improved Boole's case, this parameter was set to $\epsilon_i = \frac{\epsilon+P(A_i \cap ... \cap A_n)(n-1)}{n}$, where $n = 16$ PV systems. One potential direction for future work is addressing how to include this term as a variable to determine the optimal $\epsilon_i$ for each single chance constraint while still satisfying the new Boole's bound.

\begin{figure}[t!]
\begin{center}
\includegraphics[width=9.5cm ]{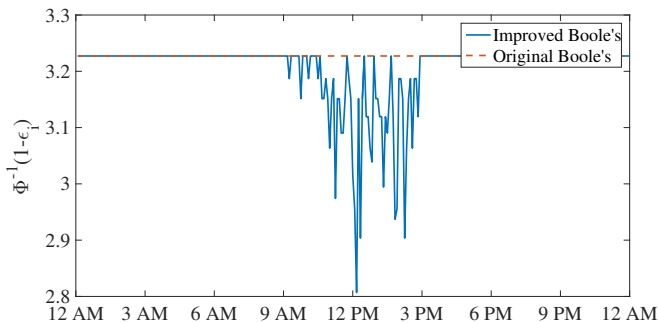}
\end{center}
\vspace{-.3cm}
\caption{Value of the inverse CDF term for August 1 under both Boole's inequality and the Improved Boole's inequality.}
\label{fig:phi}
\vspace{-.2cm}
\end{figure}

In Figure \ref{fig:phi}, the value of each of the quantile functions $\Phi^{-1}(1-\epsilon_i)$'s is shown. By using Boole's inequality, the probability of constraint violation is held constant throughout the day, and is more conservative than the case where this parameter is calculated by using the improved Boole's inequality. When the improved inequality is used to calculate $\epsilon_i$, the probability of violation is relaxed during times of solar irradiance. Because there is a nonzero probability of all considered voltages being over their limits during these times, the conservative nature of Boole's inequality can be reduced by subtracting out the intersection of these events.

\subsubsection{Cost vs. Violation tradeoff}
\begin{figure}[t!]
\begin{center}
\includegraphics[width=9.5cm ]{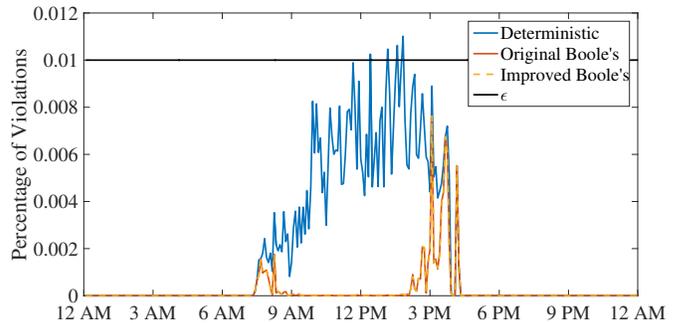}
\end{center}
\vspace{-.3cm}
\caption{Percentage of voltage violations for each Monte Carlo simulation during August 1 for each of the three cases.}
\label{fig:per_violation}
\vspace{-.2cm}
\end{figure}

In order to validate the performance of the deterministic case (e.g., using the forecasted value for PV), the case where each $\epsilon_n$ was chosen according to Boole's inequality, and the case with the improved inequality, 10,000 Monte Carlo simulations were performed for August 1-5 for each of the methods. In Figure \ref{fig:per_violation}, the total percentage of nodes with a voltage violation are shown for each of the three methods for August 1st. The deterministic case, which only considers the forecast of the PV (the mean of the random variables), violates the desired chance constraint bound of $0.01$ when compared with the chance constrained methods, because that method offers no guarantee that the voltages will be within limits with any probability. The difference between the number of voltage violations between the original Boole's method and the improved method are negligible; both bounds offer similar performance in terms of chance constraint violation probabilities.

However, the benefit of the proposed method is demonstrated in terms of cost reduction. In Table \ref{table:cost}, the total cost over the considered five-day simulation period for each of the method is shown. While the deterministic method results in the lowest cost, it also results in the highest number of violations, as well as offering little or no performance guarantee. The results from the optimization using the improved bound to calculate each $\epsilon_i$ result in a lower cost than that achieved the traditional Boole's inequality, with the cost of a small additional computation.

\begin{table}[t]
\centering
\caption{Total cost of the three methods over the five-day simulation period.}
\label{table:cost}
\begin{tabular}{l|l|l|l|l|l|l|}
\cline{2-7}
Cost(\$)                               & 8/1   & 8/2   & 8/3   & 8/4   & 8/5   & \textbf{Total} \\ \hline
\multicolumn{1}{|l|}{Deterministic}    & 18.03 & 15.56 & 13.76 & 13.99 & 11.76 & \textbf{73.10} \\ \hline
\multicolumn{1}{|l|}{Original Boole's} & 22.92 & 19.95 & 17.64 & 17.86 & 15.04 & \textbf{93.41} \\ \hline
\multicolumn{1}{|l|}{Improved Boole's} & 22.81 & 19.82 & 17.48 & 17.74 & 14.97 & \textbf{92.82} \\ \hline
\end{tabular}
\end{table}

\section{Conclusion}
In this paper, we presented a method for tractable computation of joint chance constraints in probabilistic AC OPF problems that improves upon Boole's inequality, which is often used to improve the tractability of joint chance constrained problems. An alternate tighter bound formed using Fr\'{e}chet's inequality was also given for cases of large violation probabilities (for example, sensitivity studies where we wish to find the situations with high overvoltage probabilities). The improvement upon Boole's inequality hinges on the idea of subtracting the intersection of all considered events $n-1$ times, which is in general a very fast computation to make. The new bound was proven to be a valid upper bound, and in its worse case is equivalent to that of Boole's inequality. Simulation results presented here have shown that use of the new bound preserves the probabilistic guarantees of the original joint chance constraint, while reducing the cost and conservativeness of Boole's inequality. While Gaussian random variables were considered in this paper, the framework is not distribution-dependent and could also utilize a distributionally robust formulation for the resulting single chance constraints.

Future work will consider sensitivity studies to first determine which nodes in the grid are susceptible to overvoltages, potentially increasing the impact of subtracting the intersection term. In addition, methodologies to consider $\epsilon_i$ as an optimization variable could further increase the effectiveness of the proposed new bound. Efficient methods of calculating the other intersection/union terms due to the inclusion-exclusion principle that Boole's inequality unnecessarily includes will also be explored.

\bibliographystyle{IEEEtran}
\bibliography{biblio.bib}

\begin{thebibliography}{10}
\providecommand{\url}[1]{#1}
\csname url@samestyle\endcsname
\providecommand{\newblock}{\relax}
\providecommand{\bibinfo}[2]{#2}
\providecommand{\BIBentrySTDinterwordspacing}{\spaceskip=0pt\relax}
\providecommand{\BIBentryALTinterwordstretchfactor}{4}
\providecommand{\BIBentryALTinterwordspacing}{\spaceskip=\fontdimen2\font plus
\BIBentryALTinterwordstretchfactor\fontdimen3\font minus
  \fontdimen4\font\relax}
\providecommand{\BIBforeignlanguage}[2]{{%
\expandafter\ifx\csname l@#1\endcsname\relax
\typeout{** WARNING: IEEEtran.bst: No hyphenation pattern has been}%
\typeout{** loaded for the language `#1'. Using the pattern for}%
\typeout{** the default language instead.}%
\else
\language=\csname l@#1\endcsname
\fi
#2}}
\providecommand{\BIBdecl}{\relax}
\BIBdecl

\bibitem{Roald13}
L.~Roald, F.~Oldewurtel, T.~Krause, and G.~Andersson, ``Analytical
  reformulation of security constrained optimal power flow with probabilistic
  constraints,'' in \emph{IEEE PowerTech Conference}, Grenoble, France 2013.

\bibitem{Bienstock14}
D.~Bienstock, M.~Chertkov, and S.~Harnett, ``Chance-constrained optimal power
  flow: Risk-aware network control under uncertainty,'' \emph{{SIAM} Review},
  vol.~56, no.~3, pp. 461--495, 2014.

\bibitem{Zhang_CC2011}
H.~Zhang and P.~Li, ``Chance constrained programming for optimal power flow
  under uncertainty,'' \emph{IEEE Transactions on Power Systems}, vol.~26,
  no.~4, pp. 2417--2424, Nov 2011.

\bibitem{BakerNAPS16}
K.~Baker, E.~Dall`Anese, and T.~Summers, ``Distribution-agnostic stochastic
  optimal power flow for distribution grids,'' in \emph{IEEE North American
  Power Symposium}, Denver, CO, September 2016.

\bibitem{DallAneseCDC16}
E.~Dall`Anese, K.~Baker, and T.~Summers, ``Adaptive optimal power flow for
  distribution systems under uncertain forecasts,'' in \emph{2016 Conference on
  Decision and Control (CDC)}, Las Vegas, NV, December 2016.

\bibitem{BakerTSE16}
K.~Baker, G.~Hug, and X.~Li, ``Energy storage sizing taking into account
  forecast uncertainties and receding horizon operation,'' \emph{IEEE Trans. on
  Sustainable Energy}, vol.~PP, no.~99, Aug 2016.

\bibitem{Vrak12}
M.~Vrakopoulou, K.~Margellos, J.~Lygeros, and G.~Andersson, ``Probabilistic
  guarantees for the n-1 security of systems with wind power generation,'' in
  \emph{PMAPS 2012}, Istanbul, Turkey 2012.

\bibitem{Hojjat15}
M.~Hojjat and M.~H. Javidi, ``Chance-constrained programming approach to
  stochastic congestion management considering system uncertainties,''
  \emph{IET Generation, Transmission Distribution}, vol.~9, no.~12, pp.
  1421--1429, 2015.

\bibitem{Boole1854}
G.~Boole, \emph{An Investigation of the Laws of Thought on Which are Founded
  the Mathematical Theories of Logic and Probabilities (1854)}.\hskip 1em plus
  0.5em minus 0.4em\relax New York, NY: Dover Publications, 1958.

\bibitem{Grosso14}
J.~M. Grosso, C.~Ocampo-Martinez, V.~Puig, and B.~Joseph, ``Chance-constrained
  model predictive control for drinking water networks,'' \emph{Journal of
  Process Control}, vol.~24, no.~5, pp. 504--516, 2014.

\bibitem{Nemirovski07}
A.~Nemirovski and A.~Shapiro, ``Convex approximations of chance constrained
  programs,'' \emph{SIAM J. on Optimization}, vol.~17, no.~4, pp. 969--996,
  2007.

\bibitem{Blackmore09}
L.~Blackmore and M.~Ono, ``Convex chance constrained predictive control without
  sampling,'' in \emph{Proceedings of the AIAA Guidance, Navigation and Control
  Conference}, April 2009.

\bibitem{Hong_JCC}
L.~Hong, Y.~Yang, and L.~Zhang, ``Sequential convex approximations to joint
  chance constrained programs: A monte carlo approach,'' \emph{IET Generation,
  Transmission Distribution}, vol.~59, pp. 617--630, 2011.

\bibitem{kerstingbook}
W.~H. Kersting, \emph{Distribution System Modeling and Analysis}.\hskip 1em
  plus 0.5em minus 0.4em\relax 2nd ed., Boca Raton, {FL}: {CRC} Press, 2007.

\bibitem{Liu08}
Y.~Liu, J.~Bebic, B.~Kroposki, J.~{de Bedout}, and W.~Ren, ``Distribution
  system voltage performance analysis for high-penetration {PV},'' in
  \emph{{IEEE} Energy 2030 Conf.}, Nov. 2008.

\bibitem{Tyler15}
T.~Summers, J.~Warrington, M.~Morari, and J.~Lygeros, ``Stochastic optimal
  power flow based on conditional value at risk and distributional
  robustness,'' \emph{International Journal of Electrical Power \& Energy
  Systems}, vol.~72, pp. 116--125, Nov. 2015.

\bibitem{LavaeiLow}
J.~Lavaei and S.~H. Low, ``Zero duality gap in optimal power flow problem,''
  \emph{IEEE Trans. Power Syst.}, vol.~27, no.~1, pp. 92--107, Feb. 2012.

\bibitem{PaudyalyISGT}
S.~Paudyaly, C.~A. Canizares, and K.~Bhattacharya, ``Three-phase distribution
  {OPF} in smart grids: Optimality versus computational burden,'' in \emph{2nd
  IEEE PES Intl. Conf. and Exhibition on Innovative Smart Grid Technologies},
  Manchester, UK, Dec. 2011.

\bibitem{zhang2013relaxed}
H.~Zhang, V.~Vittal, G.~Heydt, and J.~Quintero, ``A relaxed {AC} optimal power
  flow model based on a taylor series,'' in \emph{{IEEE} Innovative Smart Grid
  Technologies-Asia}, 2013.

\bibitem{madani2014}
R.~Madani, M.~Ashraphijuo, and J.~Lavaei, ``Promises of conic relaxation for
  contingency-constrained optimal power flow problem,'' in \emph{52nd Annual
  Allerton Conference on Communication, Control, and Computing}, 2014, pp.
  1064--1071.

\bibitem{coffrin2015qc}
C.~Coffrin, H.~Hijazi, and P.~{Van Hentenryck}, ``The {QC }relaxation:
  Theoretical and computational results on optimal power flow,'' \emph{IEEE
  Transactions on Power Systems}, vol.~31, no.~4, 2016.

\bibitem{sairaj2015linear}
S.~Dhople, S.~Guggilam, and Y.~Chen, ``Linear approximations to {AC} power flow
  in rectangular coordinates,'' \emph{Allerton Conference on Communication,
  Control, and Computing}, 2015.

\bibitem{swaroop2015linear}
S.~Guggilam, E.~Dall'Anese, Y.~Chen, S.~Dhople, and G.~B. Giannakis, ``Scalable
  optimization methods for distribution networks with high pv integration,''
  \emph{{IEEE} Transactions on Smart Grid}, 2016.

\bibitem{BoVa04}
S.~Boyd and L.~Vandenberghe, \emph{Convex Optimization}.\hskip 1em plus 0.5em
  minus 0.4em\relax Cambridge University Press, 2004.

\bibitem{campi}
M.~C. Campi, S.~Garatti, and M.~Prandini, ``The scenario approach for systems
  and control design,'' \emph{Annual Reviews in Control}, vol.~33, no.~1, pp.
  149--157, 2009.

\bibitem{IEEE37}
IEEE, ``37 node distribution test feeder,'' [Online] Available at
  \texttt{https://ewh.ieee.org/soc/pes/dsacom/testfeeders/}.

\bibitem{Bank13}
J.~Bank and J.~Hambrick, ``Development of a high resolution, real time,
  distribution-level metering system and associated visualization modeling, and
  data analysis functions,'' {N}ational Renewable Energy Laboratory, Tech. Rep.
  NREL/TP-5500-56610, May 2013.

\end{thebibliography}

\begin{IEEEbiography}[{\includegraphics[width=1in]{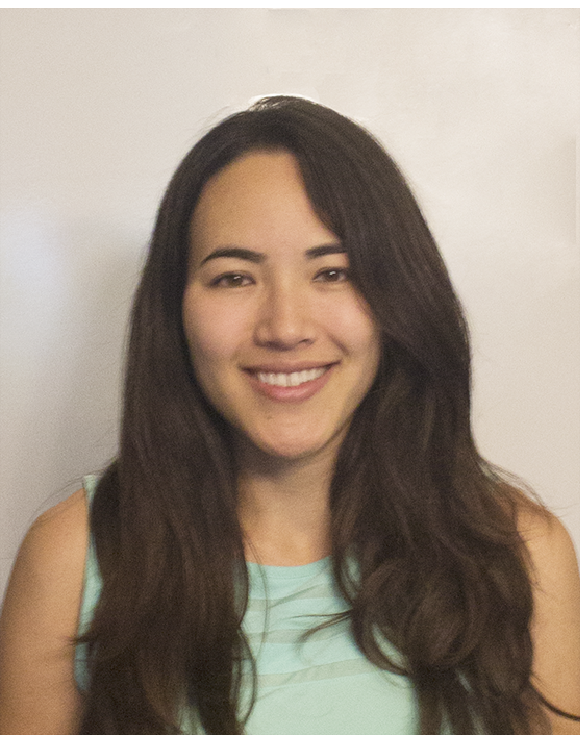}}]{Kyri Baker} (S'08, M'15)
received her B.S., M.S, and Ph.D. in Electrical and Computer Engineering at Carnegie Mellon University in
2009, 2010, and 2014. She is currently a power systems optimization and control engineer at the National Renewable Energy Laboratory in Golden, CO. Her research interests include power system optimization and planning, smart grid technologies, and renewable energy integration.
\end{IEEEbiography}

\begin{IEEEbiography}[{\includegraphics[width=1in]{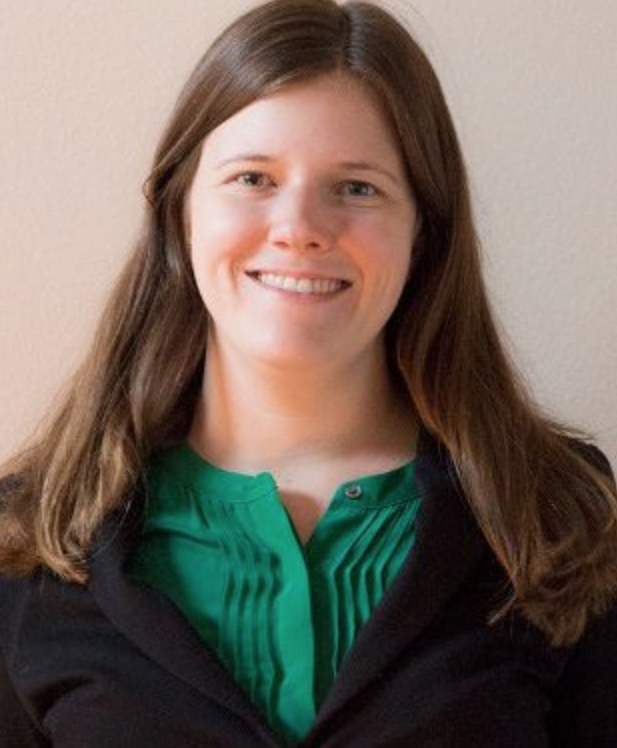}}]{Bridget Toomey}
received her B.A. in Mathematics from Grinnell College in 2013 and her M.S. in Mathematics from the University of Miami in 2015. She is currently a Content Engineer specializing in R-based predictive modeling tools at Alteryx Inc. in Broomfield, CO. Her research interests include machine learning, statistics, and data analytics in general, as well as increasing the accessibility of all three.
\vspace{-.4cm}
\end{IEEEbiography}
\vfill
\end{document}